\input ./style/arxiv-general.cfg
\documentclass[seceqn,MSNbibl,number,citesort,dvips]{arxbj}
\makeatletter
   \@ifpackageloaded{graphicx}{}{\usepackage{graphicx}}
\makeatother
\usepackage{mathrsfs}

\aid{0}
\volume{21}
\issue{4}
\pubyear{2015}
\firstpage{2120}
\lastpage{2138}
\doi{10.3150/14-BEJ637} 
\docsubty{FLA}

\makeatletter

\newtheorem{theorem}{Theorem}[section]

\newproclaim{rem}{Remark}[section]

\newproclaim{defn}{Definition}[section]

\newtheorem{cor}{Corollary}[section]

\def\R{\mathbf{R}}
\def\P{\mathbf{P}}
\def\r{\mathbf{r}}
\def\E{E}
\def\L{L}
\def\V{V}
\def\D{D}
\def\S{S}
\def\Ed{E^{\dagger}}
\def\Vd{V^{\dagger}}
\def\defeq{\stackrel{\mathrm{def}}{=}}
\def\argmin{\operatorname{arg\,min}\limits}
\newcommand{\eqref}[1]{(\ref{#1})}
\makeatother

\begin{document}
\begin{frontmatter}

\title{Maximum likelihood estimators uniformly minimize distribution variance among distribution unbiased estimators in exponential families}
\runtitle{MLEs are UMVUEs in exponential families}

\begin{aug}
\author[A]{\inits{P.}\fnms{Paul} \snm{Vos}\thanksref{e1}\ead[label=e1,mark]{vosp@ecu.edu}} \and
\author[A]{\inits{Q.}\fnms{Qiang} \snm{Wu}\corref{}\thanksref{e2}\ead[label=e2,mark]{wuq@ecu.edu}}
\address[A]{Department of Biostatistics, East Carolina University, Greenville, NC 27834, USA.\\
 \printead{e1,e2}}
\end{aug}

\received{\smonth{5} \syear{2013}}
\revised{\smonth{4} \syear{2014}}

\begin{abstract}
We employ a parameter-free distribution estimation framework where
estimators are random distributions and utilize the Kullback--Leibler
(KL) divergence as a loss function. Wu and Vos [\textit{J. Statist. Plann. Inference} \textbf{142} (2012) 1525--1536] show that when
an estimator obtained from an i.i.d. sample is viewed as a random distribution,
the KL risk of the estimator decomposes in a fashion parallel to the
mean squared error decomposition when the estimator is a real-valued
random variable. In this paper, we explore how conditional versions
of distribution expectation ($\Ed$) can be defined so that a distribution
version of the Rao--Blackwell theorem holds. We define distributional
expectation and variance ($\Vd$) that also provide a decomposition
of KL risk in exponential and mixture families. For exponential families,
we show that the maximum likelihood estimator (viewed as a random
distribution) is distribution unbiased and is the unique uniformly
minimum distribution variance unbiased (UMV$^{\dagger}$U) estimator.
Furthermore, we show that the MLE is robust against model specification
in that if the true distribution does not belong to the exponential
family, the MLE is UMV$^{\dagger}$U for the KL projection of the
true distribution onto the exponential families provided these two
distribution have the same expectation for the canonical statistic.
To allow for estimators taking values outside of the exponential family,
we include results for KL projection and define an extended projection
to accommodate the non-existence of the MLE for families having discrete
sample space. Illustrative examples are provided.
\end{abstract}

\end{frontmatter}

\section{Introduction}
\label{sec:Introduction}

Wu and Vos \cite{wv2011} introduce a parameter-free distribution estimation
framework and utilize the Kullback--Leibler (KL) divergence as a loss
function. They show that the KL risk of a distribution estimator obtained
from an i.i.d. sample decomposes in a fashion parallel to the mean squared
error decomposition for a parameter estimator, and that an estimator
is distribution unbiased, or simply unbiased, if and only if its distribution
mean is equal to the true distribution. Distribution unbiasedness
can be defined without using any parameterization. We call this approach
parameter-free even though there may be applications where it is desirable
to use a particular parameterization. When the distributions are,
in fact, parametrically indexed, distribution unbiasedness handles
multiple parameters simultaneously and is consistent under reparametrization.
Wu and Vos \cite{wv2011} also show that the MLE for distributions in the exponential
family is always distribution unbiased.

The KL expectation and variance functions $\E$ and $\V$ are defined
by minimizing over the space of all distributions. These functions
completely describe an estimator in terms of its KL divergence around
any distribution. In this paper, we introduce distribution expectation
and variance functions $\Ed$ and $\Vd$ that are defined by minimizing
over a smaller space of distributions. For exponential and mixture
families, the expected KL risk is a function only of these quantities.

Even though the focus of this paper is on parametric exponential families,
our approach is parameter-free in that the definitions and results
are provided without regard to the parameterization of the family.
There are three advantages to this approach: one, the lack of invariance
of bias across parameter transformations is avoided; two, we can allow
for estimators taking values outside of the exponential family; three,
the case where the true distribution does not belong to the family
is easily addressed.

Section~\ref{sec:Kullback-Leibler-Risk,-Variance,} introduces the
distribution expectation and variance functions and shows how these
are a generalization of the mean and expectation functions for mean
square error. Exponential families and their extension are discussed
in Section~\ref{sec:Exponential-Family}. The fundamental properties
of the distribution mean and variance functions allow using the ideas
of Rao--Blackwell \cite{blackwell1947} to show that the MLE is the
unique uniformly minimum distribution variance unbiased estimator
(UMV$^{\dagger}$UE). This result is proved in Section~\ref{sec:Rao-Blackwell-and-the}.
Three examples are given in Section~\ref{sec:examples} and Section~\ref{sec:Discussions} contains further remarks.

\section{Kullback--Leibler risk, variance, and expectation}
\label{sec:Kullback-Leibler-Risk,-Variance,}

\subsection{Motivation}

The parametric version of the Rao--Blackwell theorem can be proved
using a Pythagorean relationship that holds for mean square error (MSE) and the expectation
operator. To prove the distribution version of the Rao--Blackwell theorem,
we use a similar relationship that holds for KL risk and the KL expectation
along with a second Pythagorean relationship that holds in exponential
families for KL divergence and the KL projection. Basic properties
of the expectation operator for real-valued random variables used
in the proof can be extended to distribution-valued random variables.
We begin with the property that the expectation minimizes the MSE.

For (real-valued) random variable $Y$ and $a\in\mathbb{R,}$ we can
define the average behavior of $Y$ relative to $a$ using the risk
function
\[
\E\bigl[d (Y,a )\bigr],
\]
where $d$ is a loss function, that is, a nonnegative convex function
on $\mathbb{R}\times\mathbb{R}$. When $\E[d(Y,a)]<\infty$ for some
$a$, we define
\[
V_{d}Y\defeq\inf_{b\in\mathbb{R}}\E\bigl[d(Y,b)\bigr]
\]
and
\[
E_{d}Y\defeq\argmin_{b\in\mathbb{R}}\E\bigl[d(Y,b)\bigr]
\]
if the minimum exists, in which case,
\[
V_{d}Y=\E\bigl[d (Y,E_{d}Y )\bigr].
\]

When $d (a,b )=\L(a,b)= (a-b )^{2}$, that is, risk
is MSE, we have
%
\begin{eqnarray}
\label{eq:Eequiv}E_{L}Y&\defeq&\argmin_{b\in\mathbb{R}}E\bigl[L(Y,b)\bigr]   =   \int
y\,\mathrm{d}R_{0}\defeq EY,
\\
\label{eq:Vequiv}V_{L}Y&\defeq&\inf_{b\in\mathbb{R}}E\bigl[L(Y,b)\bigr]   =
 E
\bigl[L(Y,EY)\bigr]\defeq VY.
\end{eqnarray}
Note that we use the loss function as subscript to indicate expectation
and variance defined in terms of an argmin and infimum of the loss
function, while expectations and variances without a subscript are
defined in terms of an integral, or in terms of a sum if the sample
space is discrete. The middle equality signs in equations (\ref{eq:Eequiv})
and (\ref{eq:Vequiv}) are well-known results for $EY$ and $VY$.
These two values completely characterize the risk because of the relationship
%
\begin{equation}\label{eq:MSE1}
\E\bigl[\L (Y,a )\bigr]=\L (E_{L}Y,a )+\V_{L}Y \quad\quad\forall a
\in\mathbb{R}.
\end{equation}
In particular, the MSE for a random variable $Y$ is completely determined
by knowing its expectation $\E Y$ and variance $\V Y$. Note that
(\ref{eq:MSE1}) holds for any distribution function such that $\E Y$
and $\V Y$ exist. For general loss functions $d$, the argmin $E_{d}Y$
and min $V_{d}Y$ do not characterize the risk; that is,
\[
\E\bigl[d (Y,a )\bigr]-d (E_{d}Y,a )
\]
will be a function of $a$.

The expectation and variance also have the following conditional properties
%
\begin{eqnarray}
\label{eq:MSE2}\E Y & = & \E\E [Y|X ],
\\
\label{eq:MSE3}\V Y & = & \V\E [Y|X ]+\E\bigl[\V(Y|X)\bigr].
\end{eqnarray}
In the next section, we consider random variables that take values
on a space of distributions $\mathcal{R}$ and show that when the
KL divergence is used to compare distributions, equations (\ref{eq:Eequiv})
through (\ref{eq:MSE3}) hold for KL risk.

\subsection{Space of all distributions $\mathcal{R}$}

Let $(\mathbb{X},\mathscr{X})$ be a sample space equipped with a
$\sigma$-finite measure $\lambda$. When $\mathbb{X}$ is finite
or countable, $\lambda$ is usually the counting measure. When $\mathbb{X}\subset\mathbb{R}^{d}$
and $\mathbb{X}$ contains an open set of $\mathbb{R}^{d}$ for some
$d=1,2,\ldots $\,, then $\lambda$ is usually the Lebesgue measure on
$\mathbb{R}^{d}$. Requiring $\mathbb{X}$ to contain an open set
implies that the dimension of $\mathbb{X}$ is $d$. Let $\mathcal{R}$
be the collection of all probability measures $R$ on $(\mathbb{X},\mathscr{X})$
that are absolutely continuous with respect to $\lambda$, that is, $\lambda(A)=0$
implies $R(A)=0$ for all $A\in\mathscr{X}$. This is denoted as $R\ll\lambda$.
Note that we allow the support of $R$ to be a proper subset of $\mathbb{X}$.

Let $\mathbf{R}$ (in bold font) be a random quantity whose values
are distributions in $\mathcal{R}$. The density of the distribution
$R$ with respect to $\lambda$ will be denoted by $r$ (in lower
case), and the corresponding random variable by $\mathbf{r}$
(in bold font lower case). Following Definition 2.1 in   \cite{wv2011},
$\mathbf{R}$ is an $\mathcal{R}$-valued random variable if $\mathbf{R}(A)$
is a real-valued random variable for all $A\in\mathscr{X}$. We are
considering the problem of estimating a distribution so for this paper
$\mathbf{R}=\widehat{\mathbf{R}}_{\mathbf{X}}$ is any estimator
of an unknown distribution $R_{0}\in\mathcal{R}$ where $\mathbf{X}$
is an i.i.d. sample from $R_{0}$. A random distribution is a mapping
from $\mathbb{X}^{n}$ to $\mathcal{R}$. Let $S$ be another random
quantity that is jointly distributed with $\mathbf{R}$.

\begin{theorem}
\label{the:ptilde}For every $S=s$, $K_{s}=\E[\R|S=s]$ is a probability
measure that is absolutely continuous with respect to $\lambda$,
that is, $K_{s}\in\mathcal{R}$, is unique up to measure zero ($\lambda$),
and has a density
%
\begin{equation}\label{eq:ptilde}
k_{\mathbf{s}}(y)=\E\bigl[\mathbf{r}(y)|S=s\bigr]  \quad\quad \mbox{for }   y\in
\mathbb{X}.
\end{equation}
In addition, when $s$ is replaced with the random variable $S$,
$K_{S}=\E[\mathbf{R}|S]$ is an $\mathcal{R}$-valued random variable.
\end{theorem}

\begin{pf}
For all $s$ it is easily seen that $K_{s}$ is a probability measure
because $K_{s}$ is countably additive and $K_{s}(\mathbb{X})=1-K_{s}(\varnothing)=1$,
where $\varnothing$ is the empty set. The remaining claims of the theorem
can be established by noting that equation (\ref{eq:ptilde}) can
be written as
%
\begin{equation}\label{eq:ksintegral}
k_{s}(y)=\int_{\mathbb{X}^{n}}r_{\mathbf{x}}(y)r_{0}^{n}(
\mathbf{x}|s)\,\mathrm{d}\lambda^{n}(\mathbf{x}),
\end{equation}
where $r_{0}^{n}(\mathbf{x}|s)$ is the conditional distribution of
$\mathbf{x}$ given $s$. Since
\[
E\bigl[\mathbf{R}(A)|s\bigr]=\int_{\mathbb{X}^{n}}\int
_{A}r_{\mathbf{x}}(y)\,\mathrm{d}\lambda(y)r_{0}^{n}(
\mathbf{x}|s)\,\mathrm{d}\lambda^{n}(\mathbf{x}),
\]
the set $A\in\mathcal{X}$ is arbitrary, and the integrals can be
interchanged, we see that $k_{s}(y)$ is the density for $K_{s}$
and $K_{s}\in\mathcal{R}$ for each $s$ so $K_{S}$ is an $\mathcal{R}$-valued
random variable.
\end{pf}

For $\mathcal{R}$-valued random variable $\R$ and $R\in\mathcal{R,}$
we can define the average behavior of $\R$ relative to $R$ using
the risk function
\[
\E\bigl[d (\R,R )\bigr],
\]
where $d$ is a loss function, that is, a nonnegative convex function
on $\mathcal{R}\times\mathcal{R}$. Note that the expectation used
to define the risk is with respect to some distribution $R_{0}\in\mathcal{R}$;
$R_{0}$ will be fixed but arbitrary other than constraints to ensure
that the quantities in the expressions below exist and that the support
of $R_{0}$ is $\mathbb{X}$. For any function $d$ such that $\E[d(\R,R)]<\infty$
for some $R$, we define
\[
V_{d}\R\defeq\inf_{R_{1}\in\mathcal{R}}\E\bigl[d(
\R,R_{1})\bigr]
\]
and
\[
E_{d}\R\defeq\argmin_{R_{1}\in\mathcal{R}}\E\bigl[d(\R,R_{1})
\bigr]
\]
if the minimum exists, in which case,
\[
V_{d}\R=\E\bigl[d (\R,E_{d}\R )\bigr].
\]
For KL risk, that is, when $d(R_{1},R_{2})=D(R_{1},R_{2})\defeq E_{R_{1}}\log(r_{1}/r_{2})$,
we have
%
\begin{eqnarray}
\label{eq:Eequiv-1}E_{D}\R&\defeq&\argmin_{R_{1}\in\mathcal{R}}\E\bigl[D(\R,R_{1})
\bigr]   =   \int\r_{\mathbf{x}}(y)r_{0}^{n}(\mathbf{x})\,
\mathrm{d}\lambda^{n}(\mathbf{x})\defeq E\R,
\\
\label{eq:Vequiv-1}V_{D}\R&\defeq&\inf_{R_{1}\in\mathcal{R}}E\bigl[D(
\R,R_{1})\bigr]   =   ED(\R,E\R)\defeq V\R.
\end{eqnarray}
The middle equalities in equations (\ref{eq:Eequiv-1}) and (\ref{eq:Vequiv-1})
are established in Wu and Vos \cite{wv2011}. Since these are equal when $D$
is the KL divergence and we consider no other divergence functions
on $\mathcal{R}\times\mathcal{R}$, we will simply write $E\R\in\mathcal{R}$
and $V\R\in\mathbb{R}$ for the KL mean and variance.

Furthermore, $E\R$ and $V\R$ completely characterize the average
behavior of the $\mathcal{R}$-valued random variable $\R$ relative
to any distribution $R\in\mathcal{R}$ because of the relationship
%
\begin{equation}\label{eq:KLR1}
\E\bigl[D (\R,R )\bigr]=D (\E\R,R )+\V\R \quad\quad\forall R\in\mathcal{R}.
\end{equation}
This means the KL risk for an $\mathcal{R}$-valued random variable
$\R$, having any distribution function, is completely determined
by knowing its argmin, $\E\R\in\mathcal{R}$, and minimum, $\V\R\ge0$.
When $R=R_{0}$, equation (\ref{eq:KLR1}) gives the decomposition
of the KL risk in terms of bias and variance. The relationship in
(\ref{eq:KLR1}) will not hold for general nonnegative convex functions
$d$. In this paper we only consider KL divergence $\D(R_{1},R_{2})$.
Furthermore, a conditional expectation on $\mathcal{R}$-valued random
variables can be defined so that the following conditional properties
hold
%
\begin{eqnarray}
\label{eq:KLR2}\E\R & = & \E\E[\R|\S],
\\
\label{eq:KLR3}V\R & = & V\E[\R|\S]+\E\bigl[V(\R|\S)\bigr],
\end{eqnarray}
where $\S$ could be $\mathcal{R}$-valued but could also be real
or other valued since values of $\S$ will only be used to generate
sub sigma fields.

\begin{theorem}[(Characterization theorem for expected KL divergence on
 $\mathcal{R}$)] Let $R_{0}\in\mathcal{R}$ have support $\mathbb{X}$
 and let $\R$
be an $\mathcal{R}$-valued random variable such that the KL mean
$\E\R$ and the KL variance $\V\R$ exist and are finite. Then for
any $R\in\mathcal{R}$ the mean divergence between $\R$ and $R$
depends only on the KL mean $E\R$ and KL variance $V\R$. Furthermore,
the KL mean and KL variance satisfy the classical conditional equalities
(\ref{eq:KLR2}) and (\ref{eq:KLR3}).
\end{theorem}

\begin{pf}
Equation (\ref{eq:KLR1}) follows from the definition of KL variance
and Theorem 5.2 in   \cite{wv2011} who show that the expected KL loss
$\E[D(\R,R)]$ from an $\mathcal{R}$-valued random variable $\mathbf{R}$
to a distribution $R\in\mathcal{R}$ decomposes as
%
\begin{equation}\label{eq:decomp1}
\E\bigl[D(\R,R)\bigr]=\E\bigl[D(\R,\E\mathbf{R})\bigr]+D(\E\R,R).
\end{equation}
Equation (\ref{eq:KLR2}) follows from the fact that the KL means
$\E\R$ and $\E[\R|\S]$ have densities with respect to $\lambda$
and the order of integration can be interchanged. The steps are the
same as those that establish $\E X=\E\E[X|Y]$ for $\mathbb{R}$-valued
random variables $X$ and $Y$. We rewrite (\ref{eq:KLR1}) as
%
\begin{equation}\label{eq:V1}
\E\bigl[D(\R,R)\bigr]-D(\E\R,R)=\V\R.
\end{equation}
Note that both expectations (with domain $\mathbb{R}$-valued random
variables and with domain $\mathcal{R}$-valued random variables)
and the variance depend on the data generation distribution $R_{0}$,
which can be any point in $\mathcal{R}$ with support $\mathbb{X}$.
If this equation holds for random sample $X_{1},\ldots,X_{n}$ then
it also applies to the conditional distribution of $X_{1},\ldots,X_{n}$
given $\S=s$
\[
E\bigl[\D(\R,R)|s\bigr]-\D\bigl(E[\R|s],R\bigr)=\V(\R|s).
\]
Substituting $S$ into the equation above and taking expectation gives
%
\begin{equation}\label{eq:V3}
\E\bigl[\D(\R,R)\bigr]-\E\bigl[\D\bigl(\E[\R|\S],R\bigr)\bigr]=\E\bigl[\V(\R|\S)
\bigr].
\end{equation}
Substituting $\E[\R|\S]$ into $\R$ in (\ref{eq:V1}) and using $\E\E[\R|\S]=\E\R$
gives
%
\begin{equation}\label{eq:V4}
\E\bigl[D\bigl(\E[\R|\S],R\bigr)\bigr]-D(\E\R,R)=\V\bigl(\E[\R|\S]
\bigr).
\end{equation}
Adding (\ref{eq:V3}) to (\ref{eq:V4}) and substituting from (\ref{eq:V1})
proves (\ref{eq:KLR3}).
\end{pf}

The random variable $\R$ is a distribution function defined on the
sample space and it will be useful to relate $\R$ to a statistic
$T$. We define $\mu_{T}(R)=\E_{R}T\in\mathbb{R}^{d}$ and when we
consider only one statistic we write $\mu(R)=\mu_{T}(R)$. The $\mathbb{R}^{d}$-valued
random variable $\mu(\R)$ describes the behavior of the $\mathcal{R}$-valued
random variable $\R$ and the mean of $\mu(\R)$ can be obtained from
the KL mean.

\begin{theorem}[(Expectation property on
$\mathcal{R}$)]\label{thm:Expectation-Property-on} For any statistic $T$ such that $\mu(\R)<\infty$
a.e., the mean of $T$ under $\E\R$ equals the mean of $\mathbb{R}^{d}$-valued
random variable $\mu(\R)$
%
\begin{equation}
\mu(\E\R)=\E\bigl[\mu(\R)\bigr].\label{eq:KLExpectationProperty}
\end{equation}
\end{theorem}

\begin{pf}
The density for $\E\R$ can be written as $\int r_{\mathbf{x}}(y)r_{0}^{n}(\mathbf{x})\,\mathrm{d}\lambda^{n}(\mathbf{x})$
so that
\begin{eqnarray*}
\mu(\E\R) & = & \int T(y)\int r_{\mathbf{x}}(y)r_{0}^{n}(
\mathbf{x})\,\mathrm{d}\lambda^{n}(\mathbf{x})\,\mathrm{d}\lambda(y)
\\
& = & \int r_{0}^{n}(\mathbf{x})\int T(y)r_{\mathbf{x}}(y)
\,\mathrm{d}\lambda(y)\,\mathrm{d}\lambda^{n}(\mathbf{x})=E\bigl[\mu(\R)
\bigr]
\end{eqnarray*}
because the order of integration can be switched.
\end{pf}

\subsection{General subspace $\mathcal{P}$}

We typically are interested in a subfamily of distributions $\mathcal{P}\subset\mathcal{R}$
and we describe a distribution in terms of the KL risk $\E[\D(\R,P)]$
for $P\in\mathcal{P}$. We add the regularity condition that the support
of each distribution in $\mathcal{P}$ is $\mathbb{X}$. Equation
(\ref{eq:KLR1}) shows that $\E\R$ and $\V\R$ give the KL risk for
any $P\in\mathcal{P}$. However, generally $\E\R\notin\mathcal{P}$
even if $\R$ takes values only in $\mathcal{P}$. We consider whether
an expectation can be defined that takes values in $\mathcal{P}$
and so that (\ref{eq:KLR1}) holds. We will define this expectation
as a minimum over $\mathcal{P}$. We define
\[
\Vd\R=\inf_{P\in\mathcal{P}}\E\bigl[\D(\R,P)\bigr]
\]
and
\[
\Ed\R=\argmin_{P\in\mathcal{P}}\E\bigl[D(\R,P)\bigr]
\]
if the minimum exists, in which case
\[
\Vd\R=\E\bigl[D \bigl(\R,\Ed\R \bigr)\bigr].
\]
Equation (\ref{eq:KLR1}) now becomes
%
\begin{equation}\label{eq:KLP0}
\E\bigl[D(\R,P)\bigr]=D\bigl(\Ed\R,P\bigr)+\Vd\R+\Delta\bigl(\E\R,\Ed\R,P\bigr) \quad\quad\forall P\in
\mathcal{P},
\end{equation}
where
%
\begin{equation}\label{eq:KLcorrection}
\Delta\bigl(\E\R,\Ed\R,P\bigr)=D(\E\R,P)-D\bigl(\E\R,\Ed\R\bigr)-D\bigl(\Ed\R,P\bigr).
\end{equation}
If $\Delta$ vanishes for all $P\in\mathcal{P}$ then the argmin
$\Ed\R$ and the min $\Vd\R$ completely characterize $\R$ in terms
of KL risk. When $\Delta$ is small these functions can be used to
approximate the KL risk of $\R$. We will show the term $\Delta$
vanishes when $\mathcal{P}$ is an exponential family. The relationship
between the expectations $E\R$ and $\Ed\R$ can be expressed by using
the KL projection onto $\mathcal{P}$
\[
\Pi R=\argmin_{P\in\mathcal{P}}D(R,P).
\]
By equation (\ref{eq:KLR1}),
%
\begin{equation}\label{eq:KLProjection}
\Ed\R=\Pi\E\R.
\end{equation}
For any $\mathcal{P,}$ we have that $\V\R\le\Vd\R$ since $\mathcal{P}\subset\mathcal{R}$.
These results are summarized in the following theorem.

\begin{theorem}
Let $R_{0}\in\mathcal{R}$ such that the support of $R_{0}$ is
$\mathbb{X}$ and let $\R$ be an $\mathcal{R}$-valued random variable
such that the distribution mean $\Ed\R$ and the distribution variance
$\Vd\R$ exist and are finite. Then for any $P\in\mathcal{P}$ the
mean divergence between $\R$ and $P$ is given by (\ref{eq:KLP0}).
The term $\Delta$ measures the extent to which the KL mean, distribution
mean, and $P$ depart from forming a dual Pythagorean triangle. The
KL variance is less than or equal to the distribution variance, $\V\R\le\Vd\R$,
and the distribution mean is the KL projection of the KL mean onto
$\mathcal{P}$, $\Ed\R=\Pi\E\R$.
\end{theorem}

Wu and Vos \cite{wv2011} show that $\Delta=0$ for all $P\in\mathcal{P}$ an
exponential family. For mixture families $\E\R=\Ed\R$. Hence, $\Delta$
vanishes when $\mathcal{P}$ is either an exponential family or mixture
family.

While we don't know how to write $\Ed$ as an integral and the expectation
property (\ref{eq:KLExpectationProperty}) does not hold for $\Ed$
in general, we show equations (\ref{eq:KLR2}) and (\ref{eq:KLR3})
hold with $\E$ replaced with $\Ed$ and $\V$ replaced with $\Vd$
when $\mathcal{P}$ is either an exponential or mixture family. Furthermore,
the expectation property will hold for $\Ed$ when $\mathcal{P}$
is an exponential family and $T$ is the canonical statistic.

\section{Exponential family $\mathcal{P}$}
\label{sec:Exponential-Family}

For a general subspace $\mathcal{P}\subset\mathcal{R}$ the distribution
mean $\Ed\R$ and distribution variance $\Vd\R$ do not characterize
$\E[\D(\R,P)]$ for $P\in\mathcal{P}$. However, when $\mathcal{P}$
is an exponential family these quantities do characterize $\E[\D(\R,P)]$
and the classical equalities relating conditional mean and variance
hold. A standard reference for exponential families is Brown \cite{Brown1986},
but the approach we take here is slightly different since our emphasis
is on the distributions without regard to any particular parameterization.
An exponential family $\mathcal{P}$ will be defined by selecting
a point $P_{0}\in\mathcal{R}$ and statistic $T(x)$ taking values
in $\mathbb{R}^{d}$. The defining property of an exponential family
is that for any $P\in\mathcal{P}$ the log of the density of $P$
with respect to $P_{0}$ is a linear combination of $T(x)$ and the
constant function. We start with some definitions and basic properties.

\subsection{Definitions and the projection property}

\begin{defn}
$\mathcal{P}$ is an \emph{exponential family on} $\mathbb{X}$ if there
exists $P_{0}\in\mathcal{R}$ such that the support of $P_{0}$ is
$\mathbb{X}$ and a function $T\dvtx \mathbb{X}\mapsto\mathbb{R}^{d}$
such that for any $P\in\mathcal{P}$
\[
\mathrm{d}P\propto \mathrm{e}^{\theta'T(x)}\,\mathrm{d}P_{0}\quad\quad\mbox{for some $\theta\in
\mathbb{R}^{d}$}.
\]
The distribution $P_{0}$ is called a \emph{base point} and $T$ is
called the \emph{canonical statistic} of $\mathcal{P}$. The \emph{canonical
parameter space} is
\[
\theta(\mathcal{P})= \bigl\{ \theta\in\mathbb{R}^{d}\dvtx \mbox{for some }P
\in\mathcal{P}, \mathrm{d}P\propto \mathrm{e}^{\theta'T(x)}\,\mathrm{d}P_{0} \bigr\} .
\]
\end{defn}

Without loss of generality, we can choose a base point $P_{0}$ such
that $P_{0}\in\mathcal{P}$. We'll refer to exponential families using
base points that belong to the family.

\begin{defn}
Let $\mathcal{P}$ be an exponential family with base point $P_{0}$,
canonical statistic $T$, and set $\Theta= \{ \theta\in\mathbb{R}^{d}\dvtx
\int \mathrm{e}^{\theta'T(x)}\,\mathrm{d}P_{0}<\infty \} $.
The \emph{cumulant function} has domain $\Theta$ and is defined as
\[
\psi(\theta)=\log\int \mathrm{e}^{\theta'T(x)}\,\mathrm{d}P_{0}.
\]
The density with respect to $P_{0}$ for any $P\in\mathcal{P}$ is
\[
\frac{\mathrm{d}P}{\mathrm{d}P_{0}}=\exp \bigl\{ \theta'T(x)-\psi(\theta) \bigr\}
\quad\quad\mbox{for some } \theta\in\theta(\mathcal{P}).
\]

The family $\mathcal{P}$ is \emph{regular} if $\theta(\mathcal{P})$
is open and $\mathcal{P}$ is \emph{full }if $\theta(\mathcal{P})=\Theta$.
\end{defn}

By the factorization theorem, $T$ is sufficient. It will often be
useful to restrict the choice of $T$ so that it is complete for the
full exponential family $\mathcal{P}$.

\begin{defn}
A statistic $T$ is \emph{complete} for $\mathcal{P}$ if
\[
E_{P}h(T)=0 \quad\quad\forall P\in\mathcal{P} \quad\Longrightarrow\quad h(T)=0\quad\quad \mbox{a.e. }
\mathcal{P}.
\]
\end{defn}

The following theorem shows that the projection operator on $\mathcal{P}$
behaves like the expectation operator on $\mathcal{R}$ (Theorem \ref{thm:Expectation-Property-on})
and will be used to show that the classical conditional expectation
equation holds for $\Ed$.

\begin{theorem}[(Projection property on $\mathcal{P}$)]\label{thm:EXPmu}
If $\Pi$ is the KL projection onto $\mathcal{P}$, where $\mathcal{P}$
is an exponential family having canonical statistic $T$ and $\mu(R)=\E_{R}T$,
then for any $R\in\mathcal{R}$ such that $\mu (R )\in\mu (\mathcal{P} )$,
%
\begin{equation}\label{eq:KLproj}
\mu(\Pi R)=\mu(R),
\end{equation}
where $\mu(\mathcal{P})= \{ \mu\in\mathbb{R}^{d}\dvtx \mbox{for some }P\in\mathcal{P}, \mu=\E_{P}T \} $
is the mean parameter space of $\mathcal{P}$.
\end{theorem}

\begin{pf}
This result follows from the relationship between the natural and
expectation parameters for an exponential family $\mathcal{P}$. Let
$\mu_{1}=\mu(P_{1})$ for some $P_{1}\in\mathcal{P}$. Then the natural
parameter $\theta(P_{1})$ of this distribution satisfies
%
\begin{equation}\label{eq:dual1}
\theta(P_{1})=\operatorname{arg\,max}\limits_{\theta\in\Theta}\bigl[
\theta'\mu_{1}-\psi(\theta)\bigr]
\end{equation}
and since $\theta$ parameterizes $\mathcal{P}$,
%
\begin{equation}\label{eq:dual1-1}
P_{1}=\operatorname{arg\,max}\limits_{P\in\mathcal{P}}\bigl[\theta(P)'
\mu_{1}-\psi\bigl(\theta(P)\bigr)\bigr].
\end{equation}
The result now follows for exponential family $\mathcal{P}$ by simple
calculation
\begin{eqnarray*}
\Pi R_{1} & = & \argmin_{P\in\mathcal{P}}D(R_{1,}P)
\\
& = & \argmin_{P\in\mathcal{P}}(E_{R_{1}}\log r_{1}-E_{R_{1}}
\log p)
\\
& = & \argmin_{P\in\mathcal{P}}E_{R_{1}}\log p
\\
& = & \argmin_{P\in\mathcal{P}}\bigl(\theta(P)'\mu(R_{1})-
\psi\bigl(\theta(P)\bigr)\bigr)
\\
& = & P_{1},
\end{eqnarray*}
where $\mu(P_{1})=\mu(R_{1})$ by (\ref{eq:dual1-1}).
\end{pf}

\begin{cor}[(Pythagorean property
for exponential families)]\label{cor:Pythagorean-Property-for}
Let $\mathcal{P}$ be an exponential family
and let $R\in\mathcal{R}$ such that $\Pi R$ exists. For all $P\in\mathcal{P}$
%
\begin{equation}
D(R,P)=D(R,\Pi R)+D(\Pi R,P).\label{eq:PythagoreanExp}
\end{equation}
\end{cor}

This is a well-known result. See, for example,   \cite{cencov1982} or
  \cite{Csiszar1975}.

 We define an extended projection $\overline{\Pi}R$ to be any distribution
in $\mathcal{R}$ such that expectation and Pythagorean properties
hold and it belongs to the ``boundary'' of $\mathcal{P}$; that
is,
%
\begin{eqnarray}
\label{eq:mupi}\mu(R) & = & \mu(\overline{\Pi}R),
\\
\label{eq:extended Pythagorean}D(R,P) & = & D(R,\overline{\Pi}R)+
D(\overline{\Pi}R,P)\quad\quad \forall P\in
\mathcal{P},
\\
\inf_{P\in\mathcal{P}}D(\overline{\Pi}R,P) & = & 0.
\nonumber
\end{eqnarray}
Note that $\Pi R$ satisfies these three equalities, and that the
last two equalities imply
\[
D(R,\overline{\Pi}R)=\inf_{P\in\mathcal{P}}D(R,P).
\]
The extended projection allows us to define the extended MLE in the
next section.

The Pythagorean property allows us to improve $\mathcal{R}$-valued
random variables by the projection $\Pi$ or, more generally, by $\overline{\Pi}$.

\begin{cor}[(Projection property for $\mathcal{R}$-valued random variables)]
 If $\overline{\Pi}\R$ exists a.e., then
\[
E\bigl[D(\R,P)\bigr]\ge E\bigl[D(\overline{\Pi}\R,P)\bigr]
\]
with equality holding if and only if $\overline{\Pi}\R=\R$ a.e.
\end{cor}

\begin{pf}
Replacing $R$ with $\R$ in equation (\ref{eq:extended Pythagorean})
and taking expectations shows
\[
E\bigl[D(\R,P)\bigr]=E\bigl[D(\R,\overline{\Pi}\R)\bigr]+E\bigl[D(\overline{\Pi}
\R,P)\bigr] \quad\quad\forall P\in\mathcal{P}
\]
and the result follows from the fact that $E[D(\R,\overline{\Pi}R)]\ge0$
with equality holding if and only if $\R=\overline{\Pi}\R$ a.e.
\end{pf}

\subsection{Fundamental equations for distribution mean and variance}

For exponential families, the distribution expectation and variance
have the same properties as the KL expectation and variance. One distinction
is that the expectation property of $\E$ holds for any statistic
while for $\Ed$ the expectation property holds only for the canonical
statistic $T$.

\begin{theorem}[(Characterization of expected KL divergence on  $\mathcal{P}$)]
Let $R_{0}\in\mathcal{R}$ have support $\mathbb{X}$ and let $\R$
be an $\mathcal{R}$-valued random variable such that the distribution
mean $\Ed\R$ exists and the distribution variance $\Vd\R$ is finite.
Then for any $P\in\mathcal{P}$, where $\mathcal{P}$ is an exponential
family, the mean KL divergence between $\R$ and $P$ depends only
on the distribution mean and distribution variance
%
\begin{equation}\label{eq:KLP1}
\E\bigl[D(\R,P)\bigr]=D\bigl(\Ed\R,P\bigr)+\Vd\R \quad\quad\forall P\in\mathcal{P}.
\end{equation}
 Assuming the conditional expectations and variances exist, the distribution
mean and distribution variance satisfy the classical conditional equalities
%
\begin{eqnarray}
\label{eq:KLP2}\Ed\R & = & \Ed\Ed[\R|\S],
\\
\label{eq:KLP3}\Vd\R & = & \Vd\Ed[\R|\S]+\E\bigl[\Vd(\R|\S)\bigr],
\end{eqnarray}
 where $S$ is a real-valued random vector. Furthermore, the expectation
property holds for the canonical statistic $T$
%
\begin{equation}\label{eq:KLP4}
\mu\bigl(\Ed\R\bigr)=\E\bigl[\mu(\R)\bigr].
\end{equation}
\end{theorem}

\begin{pf}
By Corollary \ref{cor:Pythagorean-Property-for} and equation (\ref{eq:KLproj})
the correction term (\ref{eq:KLcorrection}) vanishes showing that
equation (\ref{eq:KLP1}) holds. Equation (\ref{eq:KLP4}) follows
from
%
\begin{equation}
\mu\bigl(\Ed\R\bigr)=\mu(\E\R)\label{eq:muEd}
\end{equation}
and the expectation property on $\mathcal{R}$ (\ref{eq:KLExpectationProperty}).
Equation (\ref{eq:muEd}) follows from the (extended) projection property
for exponential families (\ref{eq:KLproj}) and \eqref{eq:mupi} and
the relationship between $E$ and $\Ed$ (\ref{eq:KLProjection}).
Now equation (\ref{eq:KLP2}) follows from
\begin{eqnarray*}
\mu\bigl(\Ed\Ed[\R|\S]\bigr) & = & \E\bigl[\mu\bigl(\Ed[\R|\S]\bigr)\bigr]
\\
& = & \E\bigl[\mu\bigl(\E[\R|\S]\bigr)\bigr]
\\
& = & \mu\bigl(\E\E[\R|\S]\bigr)
\\
& = & \mu(\E\R)
\\
& = & \mu\bigl(\Ed\R\bigr),
\end{eqnarray*}
where the first equality follows from (\ref{eq:KLP4}), the second
and fifth equalities follow from (\ref{eq:muEd}), the third equality
follows from the expectation property of the KL mean on $\mathcal{R}$,
and the fourth equality follows from the conditional expectation property
that holds on $\mathcal{R}$ (\ref{eq:KLR2}). Equation (\ref{eq:KLP3})
follows again the same steps that justified (\ref{eq:KLR3}). We rewrite
(\ref{eq:KLP1}) as
%
\begin{equation}\label{eq:V1-1}
\E\bigl[D(\R,R)\bigr]-D\bigl(\Ed\R,R\bigr)=\Vd\R.
\end{equation}
If this equation holds for random sample $X_{1},\ldots,X_{n}$ then
it also applies to the conditional distribution of $X_{1},\ldots,X_{n}$
given $\S=s$
\[
\E\bigl[D(\R,R)|s\bigr]-D\bigl(\Ed[\R|s],R\bigr)=\Vd(\R|s).
\]
Substituting $\S$ into the equation above and taking expectation
gives
%
\begin{equation}\label{eq:V3-1}
\E\bigl[D(\R,R)\bigr]-\E\bigl[D\bigl(\Ed[\R|\S],R\bigr)\bigr]=\E\bigl[\Vd(\R|\S)
\bigr].
\end{equation}
Substituting $\Ed[\R|\S]$ into $\R$ in (\ref{eq:V1-1}) and using
$\Ed\Ed[\R|\S]=\Ed\R$ gives
%
\begin{equation}\label{eq:V4-1}
\E\bigl[\D\bigl(\Ed[\R|\S],R\bigr)\bigr]-D\bigl(\Ed\R,R\bigr)=\Vd\Ed[\R|\S].
\end{equation}
Adding (\ref{eq:V3-1}) to (\ref{eq:V4-1}) and substituting from
(\ref{eq:V1-1}) proves (\ref{eq:KLP3}).
\end{pf}

\section{Rao--Blackwell and the MLE as the unique \texorpdfstring{UMV$^{\dagger}$U}{UMV$^{dagger}$U} distribution estimator}\label{sec:Rao-Blackwell-and-the}

An immediate corollary to the characterization theorem on $\mathcal{P}$
(equations (\ref{eq:KLP1}), (\ref{eq:KLP2}), and (\ref{eq:KLP3}))
is that for any random distribution $\R$ and any statistic $\S$,
the random distribution $\Ed[\R|\S]$ will have the same distribution
mean and have distribution variance less than or equal to that of
$\R$. If $\S=T$ is sufficient then $\Ed[\R|T]$ is an estimator
and if $T$ is also complete $\Ed[\R|T]$ will have smaller variance
than $\R$ unless they are equal with probability one. This conditional
expectation is enough to establish a Rao--Blackwell result for distribution
estimators if these were restricted to~$\mathcal{P}$. However, since
we are allowing $\mathcal{R}$-valued estimators we also need to project
the distributions onto $\mathcal{P}$ using $\overline{\Pi}$.

For an exponential family $\{P(y;\tau)\}$ having mean parameter $\tau\in\mu(\mathcal{P})=M$
and discrete sample space we typically have that $\Pr(T\in M)<1$
while $\Pr(T\in\overline{M})=1$ where $\overline{M}$ is the closure
of $M$. In this case, the MLE does not always exist. However, the
characterization theorem applies to $\mathcal{R}$-valued estimators
so we can define an estimator that equals the MLE $P(y;t)$ when it
exists and as a distribution $\bar{P}(y;t)$ such that $\mu(\bar{P}(y;t))=t$
and $\inf_{P\in\mathcal{P}}D(\bar{P},P)=0$ if $t\notin M$. The
\emph{extended MLE} as distribution estimator is
\[
\widehat{P}^{*}(y;t)=\cases{ P(y;t) &\quad if $t\in M$,
\cr
\bar{P}(y;t) &\quad if  $t\notin M$.}
\]
Unbiasedness of $\widehat{\P}^{*}$ follows from the following theorem.

\begin{theorem}[(Distribution unbiased
estimators in exponential families)]\label{thm:unbiasedEstimatorsinExp}
Let $\mathcal{P}$ be an exponential
family with complete sufficient statistic $T$ and let $\R$ be a
$\mathcal{R}$-valued random variable. The estimator $\R$ is distribution
unbiased for $P_{0}=\Pi R_{0}$ if and only if $\mu(E([\R|T])=T$
a.e.
\end{theorem}

\begin{pf}
We must show $\Pi E\R=P_{0}$ for all $P_{0}\in\mathcal{P}$ if and
only if $\mu(E[\R|T])=T$ a.e. for all $P_{0}\in\mathcal{P}$. Consider
the following equivalencies each of which holds for all $P_{0}\in\mathcal{P}$:
\begin{eqnarray*}
\Pi E\R & = & P_{0}
\\
\iff\quad\mu(\Pi E\R) & = & \mu(P_{0})
\\
\iff\quad\mu(E\R) & = & \mu(P_{0})
\\
\iff\quad\mu\bigl(EE[\R|T]\bigr) & = & \mu(P_{0})
\\
\iff\quad E\bigl[\mu\bigl(E[\R|T]\bigr)\bigr] & = & \mu(P_{0})
\\
\iff\quad E\bigl[\mu\bigl(E[\R|T]\bigr)\bigr] & = & E(T).
\end{eqnarray*}
The first equivalence follows because the expectation of $T$ parameterizes
$\mathcal{P}$, the second equivalence follows from the projection
property for exponential families, the third equivalence follows from
the conditional expectation defined for the KL mean, the fourth equivalence
follows from the expectation property for the KL mean, and the fifth
equivalence follows from the definition of the function $\mu$. Clearly,
$\mu(E[\R|T])=T$ a.e. implies the last equality. Since $T$ is complete
and the last equality holds for all $P_{0}\in\mathcal{P}$, this implies
\[
\mu\bigl(E[\R|T]\bigr)=T\quad\quad \mbox{a.e.}
\]
\upqed\end{pf}

\begin{theorem}[(Optimality of the MLE for exponential
families)]\label{thm:Optimality-of-the} Let $X_{1},\ldots,X_{n}$ be i.i.d. from a distribution $R_{0}\in\mathcal{R}$
such that the support of $R_{0}$ is $\mathbb{X}$. Let $\mathcal{P}$
be an exponential family with complete sufficient statistic $T$ such
that $\mu(R_{0})\in\mu(\mathcal{P})$. If $\widehat{\P}$ is the MLE
or an extended MLE that exists a.e., then $\widehat{\P}$ is distribution
unbiased for the $\Pi R_{0}$ and it is the unique uniformly minimum
distribution variance estimator among all $\mathcal{R}$-valued estimators
that are distribution unbiased for $\Pi R_{0}$ and for which the
extended projection $\overline{\Pi}\R$ exists a.e. \end{theorem}

\begin{pf}
Uniqueness and uniform minimum distribution variance follow from the
projection property for $\mathcal{R}$-valued random variables, the
characterization theorem on $\mathcal{P}$ described above, and the
unbiasedness from Theorem \ref{thm:unbiasedEstimatorsinExp}.
\end{pf}

\section{Examples}
\label{sec:examples}

\subsection{Binomial distribution}

We consider the number of events or ``successes'' in $n$ trials.
The sample space is
\[
\mathbb{X}= \{ 0,1,2,\ldots,n \} .
\]
Under the assumptions that these trials are independent and each trial
has the same success probability $0<\theta<1$, the distribution of
$X$ belongs to the $n$-binomial family
\[
\mathcal{P}= \bigl\{ P\in\mathcal{R}\dvtx P(x)=P_{\theta}(x)={n \choose x}
\theta^{x}(1-\theta)^{n-x}\mbox{ for some }0<\theta<1 \bigr\} .
\]
The MLE for the parameter $\theta$ is $\hat{\theta}=x/n$ for $x\notin\{0,n\}$
but is undefined otherwise. The extended MLE (it will correspond in
a natural way to the extended MLE distribution estimator) is $\hat{\theta}=x/n$
for all $x\in\mathbb{X}$ and it is unbiased for $\theta$. However,
it is not unbiased for other parameterizations such as the odds $\nu=\theta/(1-\theta)$,
or the log odds $\gamma=\log\nu$. When viewed as a distribution,
that is, $P_{\hat{\theta}}(x)$, equivalently, $P_{\hat{\nu}}(x)$ or
$P_{\hat{\gamma}}(x)$ (where we allow the odds $\nu$ and log odds
$\gamma$ to take values in the extended reals), the MLE is the unique
uniformly minimum distribution variance unbiased estimator. As is
common practice, we have used the same notation $\hat{\theta}$ for
both the MLE and the extended MLE.

Estimators, whether real-valued or distribution-valued, are functions
with domain $\mathbb{X}$. For the $n$-binomial family an estimator
is given by a sequence of $n+1$ values, real numbers for $\hat{\theta}$
and probability distributions for $P_{\hat{\theta}}$. For $\hat{\theta}$,
we have the sequence
%
\begin{equation}\label{eq:mleReal}
\frac{0}{n},\frac{1}{n},\frac{2}{n},\ldots,\frac{n-1}{n},
\frac{n}{n}.
\end{equation}
Let $P_{\theta_{0}}$ be a distribution in $\mathcal{P}$. If probabilities
of $P_{\theta_{0}}$ are used to assign weights to the values in~(\ref{eq:mleReal}),
then the real number that is closest to the weighted values of (\ref{eq:mleReal})
is $\theta_{0}$. That is,
\[
\theta_{0}=\argmin_{\theta\in(0,1)}E \biggl(\frac{X}{n}-\theta
\biggr)^{2}.
\]
By the Rao--Blackwell theorem, for any other sequence of $n+1$ real
numbers
%
\begin{equation}\label{eq:OtherReal}
y(0),y(1),y(2),\ldots,y(n-1),y(n)
\end{equation}
that satisfy
\[
\theta_{0}=\argmin_{\theta\in(0,1)}E \bigl(y(X)-\theta
\bigr)^{2},
\]
the realized minimum will be greater than the minimum obtained using
the values in (\ref{eq:mleReal}) unless the sequences are equal,
$y(x)=x/n$ for $x\in \{ 0,1,2,\ldots,n \} $.

A distribution estimator $P_{\hat{\theta}}$ obtained from the real
valued estimator given in (\ref{eq:mleReal}) can be defined as
%
\begin{equation}\label{mleDist}
I_{0}(x),P_{1/n}(x),P_{2/n}(x),
\ldots,P_{(n-1)/n}(x),I_{1}(x),
\end{equation}
where $I_{a}$ is the indicator function for its subscript; that is,
the degenerate distribution putting all mass on 0 or 1. Since $\inf_{P\in\mathcal{P}}D(I_{a},P)=0$
it is easily checked that $\overline{\Pi}I_{a}=I_{a}$ which means
that the sequence in (\ref{mleDist}) is the extended MLE $\widehat{\P}^{*}$.
Hence, $\widehat{\P}^{*}=P_{\hat{\theta}}$. Again, we let $P_{\theta_{0}}$
be any distribution in $\mathcal{P}$. If $P_{\theta_{0}}$ is used
to assign weights to the distributions in (\ref{mleDist}), then the
distribution in $\mathcal{P}$ that is closest to the weighted average
of the distributions in (\ref{mleDist}) is $P_{\theta_{0}}$. That
is,
\[
P_{\theta_{0}}=\argmin_{P\in\mathcal{P}}E\bigl[D(P_{\hat{\theta}},P)\bigr].
\]
By the distribution version of the Rao--Blackwell theorem (Theorem
\ref{thm:Optimality-of-the}) for any estimator $\tilde{\theta}$,
expressed as a distribution estimator,
%
\begin{equation}\label{eq:OtherDist}
P_{\tilde{\theta}(0)},P_{\tilde{\theta}(1)},\ldots,P_{\tilde{\theta}(n)}
\end{equation}
that satisfies
\[
P_{\theta_{0}}=\argmin_{P\in\mathcal{P}}E\bigl[D(P_{\tilde{\theta}},P)\bigr],
\]
the realized minimum will be greater than that of the MLE (\ref{mleDist})
unless the two sequences of functions (\ref{mleDist}) and (\ref{eq:OtherDist})
are equal. Theorem \ref{thm:Optimality-of-the} provides a stronger
result than this since the distributions need not belong to $\mathcal{P}$.
In the class of all distribution unbiased estimators of the form
\[
R_{0}(x),R_{1}(x),R_{2}(x),\ldots,R_{n-1}(x),R_{n}(x)
\]
for which the extended projections $\overline{\Pi}$ exists, the MLE
(\ref{mleDist}) has smallest distribution variance. In the Hardy--Weinberg
model estimators that do not belong to the family $\mathcal{P}$ have
been suggested. We consider the details in Section~\ref{sub:Hardy-Weinberg-Model}.

\begin{figure}[b]

\includegraphics{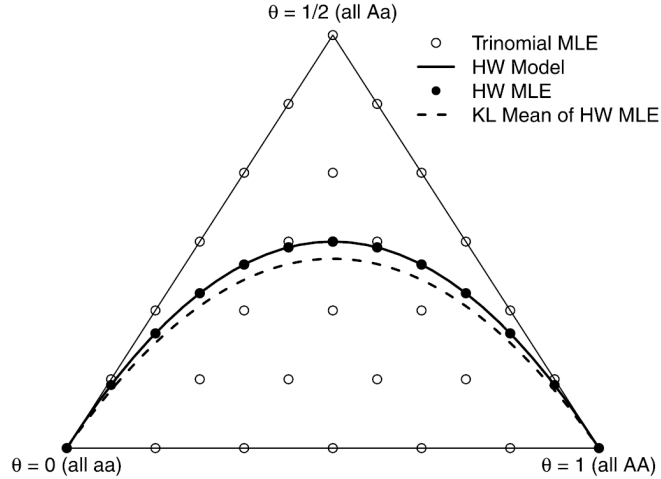}

\caption{A Hardy--Weinberg (HW) model with
$n=6$ trials. The simplex represents the trinomial model space on
$(\pi_{1},\pi_{2},\pi_{3})$ for $\pi_{1}+\pi_{2}+\pi_{3}=1$, while
the solid curve is the HW model space on $\pi_{1}(\theta)=\theta^{2}$,
$\pi_{2}(\theta)=2\theta(1-\theta)$, and $\pi_{3}(\theta)=(1-\theta)^{2}$
for $0<\theta<1$. The open circles represent the (extended) MLE under
the trinomial model $(\hat{\pi}_{1},\hat{\pi}_{2},\hat{\pi}_{3})=(Y_{1},Y_{2},Y_{3})/n$,
and the solid dots are the (extended) MLE under the HW model $\hat{\theta}=(2Y_{1}+Y_{2})/2n$.
The dashed curve shows the KL mean of the HW MLE for each value of
$\theta$.}\label{fig:Hardy-Weinberg-model}
\end{figure}

The choice of the $n$-binomial model $\mathcal{P}$ was based on
the assumptions that the data represented independent and identical
trials. If either of these assumptions were grossly violated, the
binomial model would not be appropriate. However, this model can be
used when these assumptions hold approximately in the sense that there
is a distribution $P_{0}=\Pi R_{0}$ in $\mathcal{P}$ that is close
to the data generation distribution $R_{0}$, that is, $D(R_{0},P_{0})$
is small. In this case, the MLE is the unique UMV$^{\dagger}$U estimator
for $P_{0}$.

\subsection{Hardy--Weinberg model}
\label{sub:Hardy-Weinberg-Model}

For a single pair of alleles A and a, which occur with probabilities
$\theta$ and $(1-\theta)$ for $\theta\in(0,1)$, the Hardy--Weinberg
(HW) model defines the relative frequency of genotypes AA, Aa, and
aa to be $\pi_{1}(\theta)=\theta^{2}$, $\pi_{2}(\theta)=2\theta(1-\theta)$,
and $\pi_{3}(\theta)=(1-\theta)^{2}$. For this example, we can take
$\mathcal{R}$ to be the collection of trinomial models with probabilities
$(\pi_{1},\pi_{2},\pi_{3})$ for $\pi_{1}+\pi_{2}+\pi_{3}=1$ which
can be represented by the simplex in 2-dimensional space. See Figure~\ref{fig:Hardy-Weinberg-model} for the simplex. The open circles
in Figure~\ref{fig:Hardy-Weinberg-model} are the extended MLE $(\hat{\pi}_{1},\hat{\pi}_{2},\hat{\pi}_{3})=(Y_{1},Y_{2},Y_{3})/n$
for the trinomial with $n=6$ trials, where $Y_{1}$ and $Y_{2}$
are the counts for AA and Aa. The solid curve in the simplex is the
HW model
\[
\mathcal{P}=\bigl\{(\pi_{1},\pi_{2},\pi_{3})\dvtx
\pi_{1}=\theta^{2},\pi_{2}=2\theta(1-\theta),
\pi_{3}=(1-\theta)^{2}\bigr\}
\]
which is a one dimensional exponential family with canonical sufficient
statistic $T=2Y_{1}+Y_{2}$ and canonical natural parameter $\log(\theta/(1-\theta))$.
Chow and Fong \cite{cf:1992} find the UMVU for $\pi_{1}$ and $\pi_{3}$ using
\[
E_{\theta}\bigl[\bigl(\hat{\pi}_{1}-\theta^{2}
\bigr)^{2}\bigr]+E_{\theta}\bigl[\bigl(\hat{\pi}_{3}-(1-
\theta)^{2}\bigr)^{2}\bigr]
\]
as squared-error loss. They show the UMVU is inadmissible by exhibiting
a dominating estimator. Both the UMVU and the dominating estimator
take values outside the HW model. In terms of distribution estimators,
these are $\mathcal{R}$-valued estimators.

The extended MLE for the HW model is $\hat{\theta}=(2Y_{1}+Y_{2})/2n$
while the extended distribution MLE is $P_{\hat{\theta}}$ where $P_{0}$
is the degenerate distribution putting all its mass on $(0,0,6)$ (the
lower left vertex) and $P_{1}$ is the degenerate distribution putting
all its mass on $(6,0,0)$ (the lower right vertex). The extended HW
MLE is represented by the solid dots in Figure~\ref{fig:Hardy-Weinberg-model}.

Among the difficulties with the UMVU estimator and the dominating
estimator is that there are other ways to define squared-error loss
(using one bin or two other bins). These are avoided by using KL divergence.
Since $\mathcal{P}$ is an exponential family the extended MLE is
the UMV$^{\dagger}$U for all $\mathcal{P}$-valued estimators but
also for all $\mathcal{R}$-valued estimators since the projection
exists for all points in the simplex other than the two lower vertices
which satisfy the extended projection. As a comparison, the KL mean,
represented by the dashed curve in Figure~\ref{fig:Hardy-Weinberg-model},
lives outside the model so the extended MLE isn't KL unbiased. This
is due to the curvature in the exponential family.

\subsection{Poisson distribution}

The Poisson family of distributions is
\[
\mathcal{P}=\biggl\{P\in\mathcal{R}\dvtx P_{\lambda}(x)=\mathrm{e}^{-\lambda}
\frac{\lambda^{x}}{x!} \mbox{ for some }\lambda>0\biggr\},
\]
where $x\in\mathbb{X}=\{0,1,2,\ldots\}$.

Let $X_{1},\ldots,X_{n}$ be a simple random sample from a Poisson
distribution $P_{\lambda_{0}}$. The sum $S_{n}=X_{1}+\cdots+X_{n}$
is a complete sufficient statistic of the family. Although the Poisson
family is typically parametrized by a single parameter, we consider
estimates for the probability $\Pr(X_{1}=i)=\lambda_{0}^{i}\mathrm{e}^{-\lambda_{0}}/i!$
for some $i=0,1,\ldots $\,. A crude but unbiased estimator is
\[
\delta_{0i}=\cases{ 1 &\quad if $X_{1}=i$,
\cr
0 &\quad otherwise. }
\]
Given the sum $S_{n}$, $X_{1}$ is distributed as a binomial($S_{n}$,
$1/n$) random variable, the Rao--Blackwell theorem shows that
\[
\delta_{1i}=E[\delta_{0i}|S_{n}]=\cases{\displaystyle
{S_{n} \choose i} \biggl(\frac{1}{n} \biggr)^{i}
\biggl(1-\frac{1}{n} \biggr)^{S_{n}-i} &\quad if $i\le S_{n} $,\vspace*{2pt}
\cr
0 &\quad otherwise, }
\]
is an unbiased estimator of $\Pr(X_{1}=i)$. Since $\delta_{1i}$
depends on the complete sufficient statistic $S_{n}$ only, it must
be the unique MVUE of $\Pr(X_{1}=i)$. Using the criterion of distribution
unbiasedness, these anomalous estimators do not arise. Since $S_{n}$
is the canonical statistic, the MLE $\bar{X}=S_{n}/n$ is the unique
UMVU estimator for $\lambda$ and the extended distribution MLE $P_{\bar{X}}$
is the UMV$^{\dagger}$U estimator for $P_{\lambda}$ where $P_{\bar{X}}$
is $I_{0}$ when $\bar{X}=0$.

To show how the UMVU estimator can fail completely, Lehmann \cite{Lehmann1983}
considers the parameter $\delta=(P(X=0))^{3}$ for $n=1$. In this
case, the unique UMVU estimator is $(-2)^{x}$. Since the sample consists
of nonnegative integers this estimator is represented by the following
sequence of real numbers
\[
1,-2,4,-8,16,\ldots.
\]
Parametric unbiasedness means that if the Poisson distribution that
assigns probability  $\delta^{1/3}$ to $P(X=0)$ is used to assign
probability to the terms in the sequence then $\delta=\operatorname{arg\,min}_{a\in\mathbb{R}}E((-2)^{X}-a)^{2}$.
That is, the parameter is the real number that is closest to this
sequence in terms of mean square error. In addition, the weighted
average of the above sequence is~$\delta$.

By focusing on distributions rather than the parameters that name
the distributions these problems are avoided. The MLE, as a distribution
estimator, is represented by the following sequence of probability
distributions
\[
I_{0}(x),\mathrm{e}^{-1}\frac{1^{x}}{x!},
\mathrm{e}^{-2}\frac{2^{x}}{x!},\mathrm{e}^{-3}
\frac{3^{x}}{x!},\ldots.
\]
Distribution unbiasedness means that if the Poisson distribution $P_{\lambda}$
is used to assign probability to the terms in the sequence then
\[
P_{\lambda}=\argmin_{P\in\mathcal{P}}\E\bigl[\D(P_{\hat{\lambda}},P)\bigr].
\]
That is, the distribution that generates the data is the distribution
in the exponential family that is closest to this sequence in terms
of KL risk. Any other sequence of distributions with this property
will have greater distribution variance.

\section{Discussion}
\label{sec:Discussions}

The distribution version of the Rao--Blackwell theorem \ref{thm:Optimality-of-the}
has been developed by analogy with important properties of mean square
error for the parametric version. In particular, we have used a Pythagorean-type
property for two asymmetric distribution-like functions: the KL divergence
$D(\cdot,\cdot)$ and its expectation $E[D(\cdot,\cdot)]$. For exponential
family $\mathcal{P,}$ we have
\[
D(\R,P)=D(\R,\Pi\R)+D(\Pi\R,P) \quad\quad\forall P\in\mathcal{P}
\]
while for all $\mathcal{R}$
\[
\E\bigl[D(\R,R)\bigr]=\E\bigl[D(\R,\E\R)\bigr]+\E\bigl[D(\E\R,R)\bigr]
\]
so that the expectation operator $E$ defined on $\mathcal{R}$-valued
random variables for the KL risk plays the role of the projection
operator $\Pi$ for the KL divergence. Each operator is a map from
a more complicated space to a simpler space, $E$ from $\mathcal{R}$-valued
random variables to a distribution in $\mathcal{R}$ and $\Pi$ from
distributions in $\mathcal{R}$ to a distribution in $\mathcal{P}$,
that preserve the KL risk and KL divergence, respectively.

The restriction to exponential families is essentially required by
the criterion of having a sufficient statistic of fixed dimension
for all sample sizes $n$. Specifically, the Darmois--Koopman--Pitman
theorem which follows from independent works of Darmois \cite{darmois},
Koopman \cite{koopman} and Pitman \cite{pitman} shows that when only continuous
distributions are considered, the family of distributions of the sample
has a sufficient statistic of dimension less than $n$ if and only
if the population distribution belong to the exponential family. Denny \cite{denny72}
shows that for a family of discrete distributions, if there is a sufficient
statistic for the sample, then either the family is an exponential
family or the sufficient statistic is equivalent to the order statistics.

The MLE is parameter-invariant which means that the same distribution
is named by the parametric ML estimate regardless of the parameter
chosen to index the family. One approach to studying parameter-invariant
quantities is to use differential geometry (e.g., Amari \cite{Amari1990}
or Kass and Vos \cite{Kass1997}). The parameter-invariant approach does not work
well for parameter-dependent quantities such as bias and variance
of parametric estimators. Our approach allows for the definition of
parameter-free versions of bias and variance. Furthermore, the distribution
version of the Rao--Blackwell provides two extensions: (1) minimum
variance is taken over a larger class of estimators that includes
estimators that are not required to take values in the model space
$\mathcal{P}$, (2) the true distribution need not belong to $\mathcal{P}$.

The fact that the MLE is the unique uniformly minimum distribution
variance unbiased estimator for exponential families distinguishes
the MLE from other estimators. This is in contrast to asymptotic methods
applied to MSE that can be used to show superior properties of the
MLE but, being asymptotic results, do not apply uniquely to the MLE.

Asymptotically, MSE and KL risk are the same and the MSE can be viewed
as an approximation to KL risk for large $n$. The distribution version
of the Rao--Blackwell Theorem \ref{thm:Optimality-of-the} provides
support for Fisher's claim of the superiority of the MLE even in small
samples.


\section*{Acknowledgements}
We thank the associate editor and external reviewers for their insightful
comments and suggestions which have made great improvement on this
paper.



\printhistory
\end{document}